\documentclass[12pt]{article}
\usepackage{amsthm}
\newtheorem{thm}{Theorem}

\newtheorem{obs}{Observation}
\newtheorem*{obs*}{Observation}

\newtheorem{conjecture}{Conjecture}
\newtheorem*{thm*}{Theorem}
\newtheorem{definition}{Definition}
\newtheorem*{def*}{Definition}
\newtheorem*{lemma*}{Lemma}
\usepackage{amssymb}
\usepackage{mathtools}
\usepackage{graphicx}
\usepackage{caption}
\usepackage{enumitem}

\def \nmr {\begin{enumerate}}
\def \enmr {\end{enumerate}}

\def \tmz {\begin{itemize}}
\def \etmz {\end{itemize}}

\newtheorem{claim}{Claim}
\newcommand{\Proof}{\noindent\textbf{Proof. }}
\newcommand{\smallqed}{{\tiny ($\Box$)}}

\newcommand{\s}{{\rm s}}
\newcommand{\sd}{{\rm S}}

\title{General upper bound\\ on the game domination number}

\begin{document}

\author{Csilla Bujt\'as\thanks{Research supported by the Sloveniain Research Agency under the project N1-0108}}
\date{\empty}
\maketitle

\vspace{-5ex}

 \begin{center}
 	\small Faculty of Mathematics and Physics, University of Ljubljana\\
 	Ljubljana, Slovenia\\
 	\texttt{csilla.bujtas@fmf.uni-lj.si}
 \end{center}

\begin{abstract}
It is conjectured that the game domination number is at most $3n/5$ for every $n$-vertex graph which does not contain isolated vertices. It was proved in the recent years that the conjecture holds for several graph classes, including the class of forests and that of graphs with minimum degree at least two. Here we prove that the slightly bigger upper bound $5n/8$ is valid for every isolate-free graph.
\end{abstract}


\section{Introduction}

The domination game and the corresponding graph invariant $\gamma_g(G)$ was introduced by Bre{\v{s}}ar, Klav{\v{z}}ar, and Rall in 2010 \cite{BKR-2010}. This game has been studied in many further papers,
see e.g.\ \cite{BDKK-2014, DKR-2015,  KWZ-2013, KR-2019, Kos-2017, NSS-2016, Sch-2016, XL-2018, XLK-2018}.
The notion also inspired the introduction of the total domination game 
\cite{bujtas-2018, DH-2016, HKR-2015, HKR-2017,  HKR-2018,  HR-2016, Ir-2019, JL-2019},  connected \cite{borowiecki-2019, BDIK-2019}, fractional \cite{BT-2019}, and disjoint \cite{BT-2016}  domination games, Z-, L-, LL-games \cite{ZL} on graphs,  transversal game \cite{BHT-2016, BHT-2017} and domination game on hypergraphs \cite{BPTV}.

In this paper, we prove that $\gamma_g(G) \le 5n/8$ holds for every isolate-free graph $G$.

\paragraph{Standard definitions.}
We consider simple undirected graphs. The vertex set and edge set of a graph $G$ are denoted by $V(G)$ and $E(G)$ respectively. For a vertex $v \in V(G)$, the \emph{closed neighborhood} $N[v]$ contains $v$ and its neighbors. For a set $S \subseteq V(G)$, the analogous notation $N[S]= \bigcup_{v \in S} N[v]$ is used. Then, the \emph{degree} of $v$ is $d(v)=|N[v]|-1$. If $d(v)=0$ then $v$ is an \emph{isolated vertex} and it is a \emph{leaf} if $d(v)=1$. The notations $\delta(G)$ and $\Delta(G)$ stand for the \emph{minimum} and \emph{maximum vertex degree} in $G$. If $\delta(G)\ge 1$ then the graph is \emph{isolate-free}. $P_n$ and $C_n$ respectively denote the path and cycle of order $n$. Note that $P_1$ corresponds to an isolated vertex.

A vertex dominates itself and its neighbors. A set $D \subseteq V(G)$ is a \emph{dominating set} if every vertex is dominated by at least one vertex from $D$. Equivalently, $D$ is a dominating set  if $N[D]=V(G)$. The minimum cardinality of a dominating set is the \emph{domination number} $\gamma(G)$ of the graph. 

\paragraph{Domination game.} The domination game, which was introduced by 	Bre{\v{s}}ar, Klav{\v{z}}ar, and Rall in \cite{BKR-2010}, is played on a graph $G$ by Dominator and Staller who alternately select (play) a vertex from $V(G)$. In the $i^{\rm th}$ move, the choice of $v_i$ is \emph{legal} if for the vertices $v_1,\ldots,v_{i-1}$ which have been played so far 
$$N[v_i] \setminus \bigcup_{j=1}^{i-1}N[v_j]\not=\emptyset$$
holds; that is if $v_i$ dominates at least one new vertex. The game ends with the $i^{\rm th}$ move $v_i$ if $\bigcup_{j=1}^{i}N[v_j]=V(G)$. In this case, we also say that $i$ is the \emph{value of the game}.  Dominator wants to minimize the value of the game, while Staller's goal is just the opposite. If Dominator starts the game and both players play optimally (according to their goals) the value of the game is the \emph{game domination number} $\gamma_g(G)$ of the graph $G$.

If Staller is the first to play in the domination game, we call it \emph{Staller-start game} and the analogous graph invariant is denoted by $\gamma_g'(G)$. It was proved in \cite{KWZ-2013} that $|\gamma_g(G) - \gamma_g'(G)| \le 1$ holds for every graph $G$.

\paragraph{3/5-conjecture.}
One of the central topics related to the domination game is the 3/5-conjecture posed by Kinnersley, West, and Zamani~\cite{KWZ-2013}.
\begin{conjecture} \label{conj}
If $G$ is an isolate-free graph of order $n$, then $\gamma_g(G) \le 3n/5$. 
\end{conjecture}
First, the conjecture was verified for forests of caterpillars \cite{KWZ-2013}. Then, the $3n/5$ upper bound was proved for the class of forests which do not contain leaves at a distance four apart \cite{bujtas-2015}. The latter result was extended for a larger subclass of forests in \cite{Sch-2016}. A bit later, Marcus and Peleg uploaded a manuscript to the arXiv \cite{MP-2016} where they propose a proof for the entire class of forests.\footnote{The proof is quite long and involved. It seems that the authors do not plan to submit it to a journal.} On the other hand, Conjecture~\ref{conj} was proved for graphs of minimum degree at least two \cite{HK-2016}.  As it turned out, graphs of minimum degree $\delta(G) \ge 3$ admit a better upper bound as $\gamma_g(G) \le 0.5574\,n$ always holds under this condition \cite{bujtas-2015-2}. Graphs that satisfy Conjecture~\ref{conj} with equality were investigated in~\cite{BKKR-2013, HL-2017}.

Efforts were also made to establish upper bounds on $\gamma_g(G)$ which are valid for every isolate-free graph. The first such general upper bound, namely $\lceil 7n/10\rceil$,  was proved in \cite{KWZ-2013} and later it was improved to $2n/3$ in \cite{bujtas-2015-2}. The main contribution of the present manuscript is Theorem~\ref{thm:main}, which states that $\gamma_g(G) \le 5n/8$ holds for every isolate-free $G$.

\paragraph{Our approach.} We prove the upper bound $5n/8$ by using a potential function argument. For each possible state of the game, we assign colors and numerical weights to the vertices such that the sum of the weights strictly decreases with each move of the domination game and equals zero when the game ends. The main goal is to prove that Dominator's greedy strategy ensures an appropriate lower bound on the average decrease of the potential function.

\paragraph{Structure of the paper.} The manuscript is organized around the proof of the general upper bound $5n/8$. In Section~\ref{sec:2}, we make some preparations by introducing terminology and stating some basic facts related to it. Then, Section~\ref{sec:3} is devoted to the proof of the main theorem. In the short concluding section, we put some remarks on the Staller-start version of the game.

\section{Preliminaries}
\label{sec:2}

In this section, we define some basic notions and state some observations which will be used in the proof of the main theorem.

\begin{definition} Given a graph $G$ and a set $D \subseteq V(G)$, the \emph{residual graph} $G^D$ is obtained by assigning colors to the vertices and deleting some edges according to the following rules:
	\begin{itemize}
		\item A vertex $v$ is \emph{white} if $v \notin N[D]$.
		\item  A vertex $v$ is \emph{blue} if $v \in N[D]$ and $N[v] \not\subseteq N[D]$.
		\item A vertex $v$ is \emph{red} if $N[v] \subseteq N[D]$.
		\item $G^D$ contains only those edges from $G$ that are incident to at least one white vertex.
		\end{itemize}
\end{definition}

If $G^D$ is fixed, we denote by $W$, $B$, and $R$ the set of white, blue, and red vertices, respectively. By definition, these are disjoint vertex sets and  moreover, $D \subseteq R$ and $W \cup B \cup R=V(G)$ hold.  It is clear that if $D$ is the set of vertices which have been played in a domination game so far then, in the residual graph $G^D$, a vertex $v$ belongs to $W$ if it is not dominated; $v \in B$ if $v$ is already dominated but can be played as it has some undominated neighbors; $v \in R$ if it is dominated and cannot be played in the continuation of the game. Since only the white vertices must be dominated in the later moves of the game, the edges inside $B\cup R$ do not effect the continuation. Therefore, $G^D$ contains all information that is needed for the next moves in the game. Remark further, that all red vertices are isolated in $G^D$, but we keep them by technical reasons.

In a residual graph $G^D$, the \emph{white-degree} $d_W(v)$ of a vertex $v$ is the number of its white neighbors. Analogously, we sometimes refer to the \emph{blue-degree} $d_B(v)=|B \cap N(v)|$ or to the white-blue-degree $d_{WB}(v)=d_W(v) +d_B(v)$ of a vertex. The maximum of white-degrees over the sets of white and blue vertices, respectively, are denoted by $\Delta_W(W)$ and $\Delta_W(B)$. For $i \ge 0$, $W_i$ stands for the set of white vertices with $d_W(v)=i$. Similarly, for $j \ge 1$, $B_j$ denotes the set of blue vertices having white-degree $j$. A \emph{blue leaf} in $G^D$ is a vertex from $B_1$.

In the proof, we split the game into four phases, the exact definitions of which will be given in Section~\ref{sec:3}. Namely, we will have Phase~1, 2, 3, and 4 (in this order) such that some of them might be skipped. In the latter case, the length of the phase is zero. Making a difference between the phases will simplify the proof as we then can prove properties that are satisfied by the residual graph at the end of a phase, and these have consequences for the structure in later phases. We will use two different potential functions in the proof: $f(G^D)$ will be defined for Phase~1 and 2, while $F(G^D)$ will be used in Phase~3 and 4. One more subtle distinction concerns the blue vertices in $G^D$. We say that a blue vertex $v$ is \emph{light blue} if it becomes blue during Phase~1, and it is \emph{dark blue} otherwise (i.e., if $v$ is white at the end of Phase~1 and becomes blue later).

\begin{obs}
	\label{obs:1} Let $G$ be a graph and $D \subseteq V(G)$. The following statements are true for the residual graph $G^D$.
	\begin{itemize}
		\item[$(i)$] A vertex $v$ of $G^D$ can be (legally) played, if and only if $v \in W \cup B$.
		\item[$(ii)$] If\/ $D \subseteq D' \subseteq V(G)$ and a vertex $v$ is red in $G^D$, it remains red in $G^{D'}$. If $v$ is  light blue in $G^D$, then it is either light blue or red in $G^{D'}$. Similarly, if $v$ is dark blue in  $G^D$, then it is either dark blue or red in $G^{D'}$.
		\item[$(iii)$] If\/ $v$ is a white vertex in $G^D$, then none of its neighbors are red and, consequently, $d_{WB}(v)$ equals the degree of $v$ in $G$. In particular, if $v\in W$ and  $d_{WB}(v)=1$ in $G^D$, then $v$ is a leaf in $G$.
		\item[$(iv)$] If\/ $v \in R$, then $v$ is an isolated vertex in $G^D$.  
		\item[$(v)$] $D$ is a dominating set of $G$ if and only if\/ $R = V(G)$ (or equivalently, $W = \emptyset$) in $G^D$.
		\end{itemize}
\end{obs}

 A component of $G^D$ which consists of one white and two blue vertices is a \emph{BWB-component}, a component with one white and one light blue vertex is a \emph{WB$^+$-component}, and a component with one white and one dark blue vertex is a \emph{WB$^-$-component}. Sometimes we use the notation $G^{D_i}$ that refers to the residual graph obtained after the $i^{\rm th}$ move $v_i$ in the game if $i \ge 1$. For $i=0$, $G ^{D_0}$ denotes the residual graph $G^\emptyset$ that is just the graph $G$ so that all of its vertices are white.

%By definition, $f(\cH_{j^*})=0$. Moreover, as every $v\in V(\cH_0)$ and  the corresponding edge $e_v$ satisfies $f(v)+f(e_v)=22$ in $\cH_0$, we have the following inequality.

\section{Proof of the upper bound}
\label{sec:3}

We prove the following theorem here:
\begin{thm} \label{thm:main}
Let $G$ be an isolate-free graph of order $n$. Then,
$$\gamma_g(G) \le \frac{5}{8}\, n.$$
\end{thm}
\medskip

\Proof  We assign the following weights to the vertices of a residual graph:

\begin{center}
	\begin{tabular}{l|c} 
	Type of  $v$	&   $f(v)$  \\
	\hline 	
	White &  $5$ \\
		Light blue & $4$\\
		Dark blue & $3$\\
	Red &  $0$
	\end{tabular}
\end{center}
The weight of $G^D$ is defined as the sum of the weights assigned to its vertices that is,
$$f(G^D)= \sum_{v\in V(G^D)}f(v).$$
This function $f(G^D)$ will be used as a potential function in Phase~1 and 2. By Observation~\ref{obs:1} (i), (ii), and (v), $f(G^D)$ strictly decreases with each move of the domination game and equals $0$ when the game is over. Starting with $G^D$ and supposing that the next move is $v$ in the game, $\s(v)$ denotes the decrease in the potential function. That is, $\s(v)= f(G^D) - f(G^{D\cup \{v\}})$. 
We assume throughout the proof that Dominator follows the strategy of playing a vertex $v$ from $G^D$ for which the decrease in the potential function is the largest. We prove that this greedy strategy ensures that, under an arbitrary strategy of Staller, the average decrease of the potential function in a move is at least $8$.

%Let $j_*$ be the smallest index such that the $j_*^{\rm th}$ move belongs to Dominator and $\Delta_W(W)\le 2$ and $\Delta_W(B) \le 3$ holds in the residual graph $G^{D_{j_*-1}}$.  

\subsection{Phase 1}
Phase~1 starts with the first move of the game if there exists a leaf in a component of order at least $3$. It finishes with the $i^{\rm th}$ move if  $i$ is the smallest even integer such that $G^{D_i}$ does not contain three consecutive white vertices, one of them being a leaf in $G$. If such an even $i$ does not exist, then the game (and the phase) finishes with Dominator's move in Phase~1. Recall that every blue vertex which arises in this phase is a light blue vertex with a weight of $4$. In particular, in each residual graph obtained after the end of Phase 1,  every $P_3$-subgraph that is incident to a leaf of $G$ contains at least one vertex which is either light blue or red. Therefore, we get the following claim.

\begin{claim} \label{cl:light-blue}
	If $G^D$ is a residual graph in Phase $i$ of the game so that $i \ge 2$ and $v$ is a white vertex with $d_{WB}(v)=1$, then either $v$ has a light blue neighbor or $v$ has a white neighbor $u$ such that each vertex in $N[u] \setminus\{v,u\}$ is light blue.
\end{claim}

\begin{claim} \label{cl:ph1}
In Phase 1, every move of Staller decreases $f(G^D)$ by at least $5$ and every move of Dominator decreases $f(G^D)$ by at least $11$.
\end{claim} 
\Proof If Staller plays a white vertex, it becomes red and $f(G^D)$ decreases by at least $5$. If she plays a light blue vertex $v$, it dominates at least one white vertex $u$. Since $v$ is recolored red and $u$ is recolored blue or red, we have $\s(v) \ge (4-0)+(5-4)=5$. Note that, in Phase $1$, every blue vertex is light blue. On the other hand, by definition of Phase $1$, Dominator can play a white vertex $v$  which has a  white neighbor $u$ that is a leaf in $G$ and, further, $v$ has another white neighbor $u'$.  Then, playing $v$ results in the following changes in the residual graph:  $v$ and $u$ become red as $N[u]\subset N[v]$; $u'$ becomes light blue or red as $u' \in N[v]$. It follows that $\s(v) \ge 2(5-0)+(5-4)=11.$  Dominator may play a vertex of different type but, as he follows a greedy strategy,    his every move $v$ in Phase~1 results in $\s(v) \ge 11$.
\smallqed

\begin{claim} \label{cl:av1}
	If Phase 1 consists of $p_1$ moves, then $f(G^D)$ decreases by at least $8p_1$ during this phase. 
\end{claim} 
\Proof If $p_1$ is even then, by Claim~\ref{cl:ph1}, the decrease is at least $11p_1/2 + 5 p_1/2= 8 p_1$. If $p_1$ is odd, then  the decrease is at least $11(p_1+1)/2 + 5 (p_1-1)/2= 8 p_1+3$.
\smallqed

\subsection{Phase 2}
After the end of Phase~1, every vertex which turns blue becomes a dark blue vertex with weight $3$, but we keep the light blue color and the higher weight of those vertices which were already blue at the end of Phase~1. If the first phase finishes with $G^D$, the next move belongs to Phase~2 if there is a vertex $v$ such that $f(G^D)- f(G^{D\cup \{v\}}) \ge 11$. (Otherwise, Phase~2 is skipped and its length is $0$.) Phase~2 finishes with the $i^{\rm th}$ move if $i$ is the smallest even integer such that there is no vertex $v$ with the property $f(G^{D_i})- f(G^{D_i\cup \{v\}}) \ge 11$.

 %That is, Phase~2 is the part of the game, after Phase~1, when Dominator can decrease $f(G^D)$ by at least $11$ in each of his consecutive moves.
 
\begin{claim} \label{cl:av2}
	If Phase 2 consists of $p_2$ moves, then $f(G^D)$ decreases by at least $8p_2$ during this phase. 
\end{claim} 
\Proof In Phase~2, by definition, every move of Dominator decreases $f(G^D)$ by at least $11$. Concerning a move of Staller, if she plays a white vertex, it becomes red and the decrease is at least $5$; if she plays a blue vertex $v$ with a white neighbor $u$, then $v$ becomes red and $u$ becomes dark blue or red. It is also true in the latter case that $\s(v) \ge (3-0)+(5-3)=5$. Then, if $p_2$ is even, the total decrease during Phase~2 is at least $11p_2/2 + 5 p_2/2= 8 p_2$. If $p_2$ is odd, then  it is at least $11(p_2+1)/2 + 5 (p_2-1)/2= 8 p_2+3$.
\smallqed
\medskip

Now, we prove two additional claims that describe some structural properties of $G^D$ at the end of Phase~2 and some properties that remain valid in the continuation of the game. $G^D[W]$ denotes the subgraph induced by the white vertices in $G^D$.
\begin{claim} \label{cl:end2}
	Let $G^D$ be the residual graph obtained at the end of Phase~2. Then, it satisfies the following properties:
	\begin{itemize}
		\item[$(i)$] $\Delta_W(W) \le 2$ and $\Delta_W(B) \le 3$;
		\item[$(ii)$] Every component of $G^D[W]$ is isomorphic to  $P_1$ or $P_2$ or to $C_k$ with $k \ge 4$;
		\item[$(iii)$] There is no edge between the vertices of $W_0$ and $B_3$.  
	\end{itemize}
\end{claim} 
\Proof Since Phase~2 finishes with $G^D$, the decrease $\s(v)$ in $f(G^D)$ is at most $10$ for every $v \in V(G)$.

(i) If there is a white vertex $v$ with $d_W(v)\ge 3$ in $G^D$ then, after playing it, $v$ becomes red and its white neighbors become dark blue or red. This gives $\s(v) \ge 5 + 3(5-3)=11$ that is a contradiction. Similarly, if there exists a blue vertex $u$ with more than three white neighbors then $\s(u)$ would be at least $ 3 + 4(5-3)=11$.

(ii) $\Delta_W(W) \le 2$ implies that each component of $G^D[W]$ is a path or a cycle. First, assume that there is a path component $P_j=v_1v_2\dots v_j$ for an integer $j \ge 3$. In $G^{D \cup\{v_2\}}$, the vertices $v_1$ and $v_2$ become red and $v_3$ becomes dark blue or red. This gives $\s(v_2)\ge 2\cdot 5 + 2=12$ that is a contradiction.  Further, if there exists a component which is a $3$-cycle, playing any vertex $u$ from it all the three vertices are recolored red and $s(u) \ge 3\cdot 5$, a contradiction again.

(iii) Suppose that the vertices $v \in W_0$ and $u \in B_3$ are adjacent and let $u_1$ and $u_2$ be the further white neighbors of $u$. Consider $G^{D \cup\{u\}}$. In this residual graph $v$ and $u$ are red, while $u_1$ and $u_2$ are dark blue or red. We may infer $\s(u)\ge 5+ 3 +2\cdot 2 =12$ that is a contradiction.
\smallqed

\begin{claim} \label{cl:later2}
	Let $G^D$ be an arbitrary residual graph obtained in Phase~3 or Phase~4. Then, it satisfies the following properties:
	\begin{itemize}
		\item[$(i)$] $\Delta_W(W) \le 2$ and $\Delta_W(B) \le 3$;
		\item[$(ii)$] If $v$ belongs to $W_0$ in $G^D$, then it has at least one neighbor from $B_1 \cup B_2$. In particular, if $v \in W_0$, then either $v$ has a blue leaf neighbor in $G^D$ or\/ $G^{D\cup\{v\}}$ contains a blue leaf.
		\end{itemize}
\end{claim} 
\Proof Let $G^{D^*}$ be the residual graph obtained at the end of Phase~2.

(i) Since no new white vertices appear during the game and a new blue vertex may arise only from a white vertex, Claim~\ref{cl:end2} (i) clearly implies the statement.

(ii) Recall first that $d(z)\ge 1$ is supposed for every vertex $z$ of $G$ and that, by Observation~\ref{obs:1} (iii), $d_{WB}(z)=d(z)$ holds in the residual graph for each white vertex $z$. As follows, each vertex $v$ from $W_0$ has at least one blue neighbor in $G^D$.  Consider the three cases according to the status of $v$ at the end of Phase~2. If $v$ belongs to $W_0$ in $G^{D^*}$ then,  by Claim~\ref{cl:end2} (iii), all neighbors of it are blue vertices of white-degree $1$ or $2$ in $G^{D^*}$. This remains true in $G^D$. Suppose now that $v$ belongs to $W_1$ and its only white neighbor is $u$ in $G^{D^*}$. Since $v \in W_0$ holds in $G^D$, here $u$ must be blue having only one white neighbor, namely $v$. Thus, $u \in B_1$ and the statement is true for $v$. Finally, suppose that $v \in W_2$ in $G^{D^*}$ that is $v$ belongs to a white cycle component of $G^{D^*}[W]$. Then, in $G^{D^*}$, both white neighbors $v_1$ and $v_2$  have white-degree $2$. Since $v \in W_0$ in $G^D$, here both $v_1$ and $v_2$  are blue and may have white-degree at most $2$. As they remain adjacent to $v$ in the new residual graph, $v$ has a neighbor from $B_1 \cup B_2$ in $G^D$. Finally, observe that  if $v \in W_0$ and it has a neighbor $z'$ from $B_2$ in $G^D$, then $z'$ becomes a blue leaf in $G^{D\cup\{v\}}$.
\smallqed

\subsection{Phase 3}

After the end of Phase~2, we use a new potential function $F(G^D)$. To define it, we first introduce some notations. Consider the residual graph $G^{D^*}$ and the corresponding $G^{D^*}[W]$ obtained at the end of Phase~2. Let us denote by $X_1,\dots ,X_\ell$ the vertex sets of the cycle components in $G^{D^*}[W]$. The cycle induced by $X_i$ will be called \emph{X-cycle} and denoted by $C(X_i)$, for $1\le i \le \ell$. Note that, by Claim~\ref{cl:end2} (i) and (ii), $W_2=X_1 \cup \cdots \cup X_\ell$, while  $W_1$ and $W_0$ contain, respectively, the vertices from the $P_2$- and $P_1$-components of $G^{D^*}[W]$.

We fix the term `X-cycle' and the notation $X_1, \dots , X_\ell$ at the end of Phase~2 and use it later, even if some (or all) vertices from $X_i$ become blue or red. We say that an X-cycle $C(X_i)$ is \emph{closed} in a residual graph $G^D$ if all of its edges  are present in $G^D$. Remark that a closed X-cycle may contain blue vertices, but cannot contain red vertices and blue leaves.  
An X-cycle $C(X_i)$ is \emph{open} in $G^D$ if there is a vertex $v \in X_i$ which is a blue leaf in a component of order at least four, and it is \emph{finished} if all the vertices are red except those that belong to BWB-components.
We say that $C(X_i)$ is finished with a move of the game if it was open before the move and finished after it. If $C(X_i)$ is a $4$-, $5$-, or $6$-cycle, it is possible that $C(X_i)$ directly turns from closed to finished. It cannot happen if $|X_i|\ge 7$.

In a residual graph $G^D$, we denote the number of open X-cycles by $x(G^D)$, and the number of BWB- and WB$^+$-components by $c_3(G^D)$ and $c_2(G^D)$, respectively. From the beginning of Phase~3, we use the following potential function:
\begin{equation*}
\begin{split}
F(G^D) & = \left(\sum_{v\in V(G^D)}f(v)\right) -x(G^D) -c_2(G^D) -3c_3(G^D)\\
   & = f(G^D)-x(G^D) -c_2(G^D) -3c_3(G^D).
  \end{split} 
     	  \end{equation*}
We also introduce the notation $\sd(v)= F(G^D)-F(G^{D \cup \{v\}})$.
Technically, for the residual graph $G^{D^*}$ obtained after the last move of Phase~2, we calculate both values to use $f(G^{D^*})$ in Phase~2 and $F(G^{D^*})$ in Phase~3. It is clear that $f(G^{D^*}) \ge F(G^{D^*})$. Remark that the number of BWB-, and WB$^+$-components decreases if and  only if a vertex is played from such a special component. With this move $F(G^D)$ decreases by at least $8$. Concerning the decrease in $x(G^D)$, we prove the following claim.

\begin{claim} \label{cl:X-cycle-1}  
	Let $G^D$ be a residual graph from Phase~3 and suppose that a vertex $v$ is played in $G^D$.
 \begin{itemize}
  \item[$(i)$]  If $v \in W$ or $v$ is a blue vertex from an X-cycle, then $x(G^D)$ may decrease by at most $1$.
  \item[$(ii)$] For $i=1,2,3$, if $v \in B_i$, then $x(G^D)$ may decrease by at most $i$.
 \end{itemize}
\end{claim}
\Proof At the end of Phase~2, an X-cycle $C(X_i)$ consists of white vertices of white-degree $2$ that cannot be adjacent to the vertices of a different X-cycle. Then, in $G^D$, each blue vertex $u \in X_i$ is adjacent to either one or two white vertices from $X_i$ and to none of the vertices outside $X_i$. A white vertex $z \in X_i$ has exactly two neighbors from the cycle $C(X_i)$ and, additionally, it may have some outer blue neighbors, but these neighbors never belong to other X-cycles. As follows, playing a vertex from $X_i$ cannot cause changes in the colors of vertices from $X_j$, if $j \neq i$. Also, if a vertex $y \in X_j$ is adjacent to a blue vertex $y'$ which is outside the X-cycles, $y'$ remains blue (and connected to $y$) after playing any $v \in X_i$. This proves part (i).

Concerning (ii), it is enough to consider a vertex $v \in B_i$ which is outside the X-cycles. Playing $v$ changes the colors of the white neighbors $u_1, \dots , u_i$ of $v$ and may change the colors of the blue neighbors of $u_1, \dots , u_i$. On the other hand, observe that a blue neighbor $z_j \in N[u_j]$ becomes red only if  $N[z_j] \subseteq \{u_1, \dots , u_i\}$ holds in $G^D$. Therefore, if an X-cycle $C(X_k)$ does not contain a vertex from $u_1, \dots, u_i$, then the colors of the vertices in $N[X_k]$ remain unchanged in $G^{D \cup \{v\}}$. We may conclude that only the X-cycles incident to  $u_1, \dots , u_i$ can be finished with the move $v$.
\smallqed
\medskip

We say that a component of $G^D$ is \emph{special} if it is of order $2$ or a BWB-component.

\begin{claim} \label{cl:ph2-leaf}
	If there is a blue leaf in a non-special component of $G^D$, then there exists a vertex $v$ such that $\sd(v) \ge 11$.
\end{claim} 
\Proof Suppose that $u$ is a blue leaf in a non-special component and adjacent to the white vertex $u'$ in $G^D$. We prove the claim by considering three cases depending on the white-degree of $u'$.

\textit{Case 1.} $d_W(u')=2$.\\
Let $v=u'$ and observe that, in $G^{D\cup \{u'\}}$, vertices $u$ and $u'$ are red  and the two white neighbors of $u'$ belong to $B \cup R$. Since a white vertex is played, Claim~\ref{cl:X-cycle-1} implies that at most one X-cycle may be finished with this move. Hence, $\sd(u')\ge 3+5+2 \cdot 2 -1 = 11$.

\textit{Case 2.} $d_W(u')=1$.\\
Let $v$ be the white neighbor of $u'$. In $G^{D\cup \{v\}}$ all the three vertices $u$, $u'$ and $v$ are red and at most one X-cycle becomes finished. Then, we have $\sd(v)\ge 3+2 \cdot 5 -1 = 12$.

\textit{Case 3.} $d_W(u')=0$.\\
Since it is a non-special component, either $u'$ is adjacent to at least three blue leaves (including $u$) and then $\sd(u')\ge 5+3 \cdot 3 -1 = 13$ holds,  or $u'$ has a neighbor $z$ from $B_2 \cup B_3$. In the latter case, if $z \in B_2$, then $\sd(z)\ge 3+5+3 +(5-3) -2=11$; if $z \in B_3$, then $\sd(z)\ge 3+5+3 +2 (5-3) -3=12$.
\smallqed 

\begin{claim} \label{cl:x-cycle}
	Suppose that the number of open X-cycles decreases by a move $v$ in the residual graph $G^D$. 
	\begin{itemize}
	\item[$(i)$] If it is Staller's move, then $\sd(v) \ge 6$;
	\item[$(ii)$] If it is Dominator's move, then $\sd(v) \ge 11$;
	\end{itemize}	  
\end{claim} 
\Proof 
(i) Assume that $C(X_i)$ is open in $G^D$ and finished in  $G^{D\cup \{v\}}$. We consider the following cases.

\textit{Case 1.} No new BWB-component arises in $X_i$.\\
Since $C(X_i)$ is open, $G^D$ contains a component $K$ of order at least $4$ which intersects $X_i$. Then, $K \cap X_i$ includes at least one white vertex and at least two blue leaves. If no new BWB-component arises and $C(X_i)$ becomes finished, all vertices from $K \cap X_i$ turn red. On the other hand, by Claim~\ref{cl:X-cycle-1}, $x(G^D)$ may decrease by at most $3$ in one move of the game. Therefore, we have $\sd (v) \ge 5+2 \cdot 3 -3= 8$.

\textit{Case 2.} A new BWB-component arises in $X_i$.\\
Assume first that $x(G^D)$ decreases by $1$ or $2$ when $C(X_i)$ becomes finished. In this case, to finish $C(X_i)$, at least one white vertex $u$ from $X_i$ truns blue or red. If $u$ becomes red, $c_3(G^D)$ increases, and $x(G^D)$ falls,  the decrease in the potential function is at least $5 + 3 -2 =6$. If $u$ becomes blue, a neighbor $v$ was played. In $G^{D \cup \{v\}}$, $v$ is red and $u$ is blue or red. A similar calculation as before shows that $\sd(v) \ge 3+ (5-3) + 3-2 = 6$. 

It remains to prove that the inequality holds if three X-cycles, say $C(X_i)$,
 $C(X_j)$, and $C(X_k)$, are finished with the move $v$. Then, by Claim~\ref{cl:X-cycle-1}, $v$ must be a vertex from $B_3$ which does not belong to any X-cycles. If at least one new BWB-component arises in $X_i \cup X_j \cup X_k$, Case~1 can be applied for an appropriate X-cycle and the estimation follows. If all the three X-cycles become finished without getting a new BWB-component, then a calculation similar to the previous one gives $\sd(v) \ge 3 + 3(5-3)+ 3 \cdot 3 -3 = 15$. This finishes the proof of (i). 
 \medskip
 
(ii) As an open X-cycle must contain a blue leaf in a non-special component, Claim~\ref{cl:ph2-leaf} and Dominator's greedy choice directly implies $\sd(v) \ge 11$.
\smallqed

\medskip

We want to prove that the average decrease in $F(G^D)$ is at least $8$ over the moves in Phase~3. But it is possible that Staller's move gives only $\sd(v) =5$ and Dominator's move gives $\sd(v) =10$ in Phase~3. Thus, the following claim is crucial to our final estimation.

\begin{claim} \label{cl:ph2-St5}
	If Staller plays a vertex $v$ so that $F(G^D)$ decreases by exactly $5$, then Dominator can select a vertex as his next move which decreases $F(G^{D\cup \{v\}})$ by at least $11$.
\end{claim} 
\Proof Suppose that Staller plays a vertex $v$ in $G^D$ such that $\sd(v)=5$. Then, $v$ does not belong to a special component of $G^D$ and further, Claim~\ref{cl:x-cycle} implies that no X-cycles are finished with this move $v$. First, we prove that  either a blue leaf remains in a non-special component of $G^{D\cup \{v\}}$ or the residual graph contains a closed X-cycle with at least one blue vertex.

  \textit{Case 1.} Vertex $v$ is white.\\
In this case, $\sd(v)=5$ implies that each neighbor of $v$ is from $B_2 \cup B_3$. That is, $v \in W_0$ in $G^D$ and it does not have a blue leaf neighbor. Then, by Claim~\ref{cl:later2} (ii), $G^{D \cup \{v\}}$ contains a blue leaf $u$ that was a neighbor of $v$ in $G^D$. Assume now that $u$ is a blue leaf in a special component $K$ of $G^{D \cup \{v\}}$. If $K$ is a BWB- or a WB$^+$-component, then $c_3$ or $c_2$ would increase with this move and $\sd(v)$ would be at least $6$. It remains to prove that $u$ cannot belong to a WB$^-$-component. Assuming this situation, the white vertex $z$ from the WB$^-$-component has only one neighbor in  $G^{D \cup \{v\}}$ which is dark blue. By Claim~\ref{cl:light-blue}, it is not possible. Hence, the claim is verified for the first case.

  \textit{Case 2.} Vertex $v$ is blue.\\
To comply with $\sd(v) =5$, $v$ must be a dark blue leaf in $G^D$ and its white neighbor $u$ cannot be from $W_0$. Suppose first that $u$ has only one white neighbor $u'$. Then, $u$ is a dark blue leaf in  $G^{D \cup \{v\}}$ and does not belong to a WW- or WB$^+$-component in  $G^{D \cup \{v\}}$. If $u$ is in a special component BWB, then $c_3$ increases by one in this move and then, $\sd(v) = 8$, a contradiction. If $u$ and $u'$ form a WB$^-$-component in  $G^{D \cup \{v\}}$ then, by applying Claim~\ref{cl:light-blue} for $u'$ in $G^D$, we get that $v$ was a light blue vertex and $\sd(v)= 4+(5-3) =6$. This contradiction proves the statement for $d_W(u)=1$. Suppose now that $d_W(u)=2$ that is $u$ is from an X-cycle $C(X_i)$ where the neighbors $u_1$ and $u_2$ are white. In $G^{D \cup \{v\}}$, vertex $u$ is blue. Hence, if $C(X_i)$ is still closed in $G^{D \cup \{v\}}$, then it is a closed X-cycle with a blue vertex. If $C(X_i)$ is open then, by definition, it contains a blue leaf in a non-special component. Since $C(X_i)$ contains adjacent white vertices, it cannot be a finished X-cycle in $G^D$ and since it does not become finished with the move $v$, it cannot be a finished X-cycle in $G^{D \cup \{v\}}$ either.

 Now, we can complete the proof of the claim. If a blue leaf remains in a non-special component of $G^{D\cup \{v\}}$, Claim~\ref{cl:ph2-leaf} directly implies the present statement. Now suppose that $G^{D\cup \{v\}}$ contains a closed X-cycle $C(X_i)$ with a blue vertex $z \in X_i$ such that $|X_i|\ge 7$. Since the cycle is closed, both neighbors of $z$ are white. Let $z_1$ and $z_2$ be these neighbors. If $z_1 \in W_0$ then, after playing $z$, vertices $z$, $z_1$ turn red, $z_2$ turns blue or red  and moreover, $C(X_i)$ becomes open. It follows that $F( G^{D\cup \{v\}})$ decreases by at least $3+5+2+1=11$. If $z_1 \in W_1$ in  $G^{D\cup \{v\}}$, then let $z_1'$ be its white neighbor. After playing $z_1'$, both $z_1$ and $z_1'$ are recolored red and the X-cycle becomes open. This gives $\sd(z_1') \ge 2\cdot 5 +1=11$ again. For smaller X-cycles, where $|X_i|=4,5,6$, it might happen that $C(X_i)$ becomes finished and not open after the moves $z$ or $z_1'$. In these cases either a new BWB-component arises or more vertices turn red than counted above. It is easy to check that the estimations on $\sd(z)$ and $\sd(z_1')$ remain valid.
\smallqed

\begin{claim} \label{cl:av3}
	If Phase 3 consists of $p_3$ moves, then $F(G^D)$ decreases by at least $8p_3$ during this phase. 
\end{claim} 
\Proof We know that every move of Staller decreases $F(G^D)$ by at least $5$ and, by definition, every move of Dominator decreases it by at least $10$. Claim~\ref{cl:ph2-St5} implies that each move $v_i$ of Staller and the next move $v_{i+1}$ of Dominator together decrease $F(G^D) $ by at least $16$. Remark that we have $\sd(v) \ge 10$ for the first move $v$ of Phase~3 and further, by Claim~\ref{cl:ph2-St5}, the last move $u$ of Staller results in $\sd(v) \ge 6$. If $p_3$ is even, we obtain that $F(G^D)$ decreases by at least $10 + 16 \cdot \frac{p_3-2}{2}+6= 8p_3$. If $p_3$ is odd, then the entire game finishes with Dominator's move in the Phase~3. For this case, the bound can be proved similarly. \smallqed 

\begin{claim} \label{cl:end3}
	Let $G^D$ be the residual graph obtained at the end of Phase~3. Then, it satisfies the following properties:
	\begin{itemize}
		\item[$(i)$] Every X-cycle is finished and, in particular, $W_2= \emptyset$;
		\item[$(ii)$] $W_1=\emptyset$; 
		\item[$(iii)$] $B_2 \cup B_3=\emptyset$;
		\item[$(iv)$] No white vertex has more than two blue neighbors.		
	\end{itemize}
\end{claim} 
\Proof Since Phase~3 finishes with $G^D$, the decrease $\sd(v)$ in $F(G^D)$ is at most $9$ for every $v \in V(G)$. 

(i) Suppose first that the X-cycle $C(X_i)$ is closed in $G^D$. We consider the following three cases.

  \textit{Case 1.} All vertices from $X_i$ are white that is, $X_i \subseteq W_2$.\\
If $C(X_i)$ is a $4$-cycle, then for any $v \in X_i$, vertex $v$ becomes red and the two white neighbors become blue in $G^{D\cup \{v\}}$. Moreover, a new BWB-component arises and $C(X_i)$ becomes finished. Consequently, $\sd(v) \ge 5+2\cdot 2 +3=12$ that is a contradiction. If the length of $C(X_i)$ is at least $5$, select an arbitrary $u \in X_i$ and observe that $C(X_i)$ is an open cycle in $G^{D\cup \{u\}}$. This implies $\sd(v) \ge 5+2\cdot 2 +1=10$, a contradiction again. 
 
 \textit{Case 2.} There exists a vertex $v\in X_i \cap W_0$.\\
Since the cycle is closed, the two neighbors of $v$, say $v_1$ and $v_2$, belong to $B_2$. Consider $G^{D\cup \{v_1\}}$ and observe that $v$ and $v_1$ become red and the other white neighbor $u$ of $v_1$ becomes blue or red. In particular, if the cycle is finished with this move, then either a new BWB-component arises or $u$ also turns red. In the former case, we infer $\sd(v_1) \ge 3 + 5 +(5-3)+3 -1= 12$, while for the latter one we conclude $\sd(v_1) \ge 3 + 2\cdot 5 -1= 12$. If $C(X_i)$ is not finished in $G^{D\cup \{v_1\}}$, it becomes open and we get  $\sd(v_1) \ge 3 + 5 +(5-3) +1= 11$. Each subcase yields a contradiction.

  \textit{Case 3.} There exist two adjacent vertices $u$ and $v$ such that $u\in X_i \cap W_1$ and  $v\in X_i \cap (W_1 \cup W_2)$.\\
In $G^{D\cup \{v\}}$, the vertices $u$ and $v$ become red. Thus, even if $C(X_i)$ becomes finished and not open with this move, $\sd(v) \ge 2\cdot 5=10$, a contradiction. 
%If the cycle is finished, then either a new BWB-component occurs, or we had a $4$-cycle with one blue and three white vertices. These cases give the contradictions $\sd(v) \ge 10+3-1$ and $\sd(v) \ge 15+3-1$.

The above cases cover all possibilities for a closed X-cycle. Hence, it is enough to prove that $C(X_i)$ is not open in $G^D$. By definition, an open X-cycle contains a blue-leaf in a non-special component and, by Claim~\ref{cl:ph2-leaf}, there exists a vertex $v$ such that $\sd(v) \ge 11$. This contradiction and  the observation that a finished cycle cannot contain two adjacent white vertices complete the proof of (i).

\medskip

(ii) By (i), we infer $W=W_0 \cup W_1$. So, if $W_1 \neq \emptyset$, we may choose two adjacent vertices, say $v$ and $u$, from it. As also follows from (i), the number of open cycles cannot decrease when $v$ is played. Moreover, $u$ and $v$ turn red.  Therefore, we have $\sd(v) \ge 2 \cdot 5=10$, a contradiction.

\medskip
(iii) By (i) and (ii), all white vertices belong to $W_0$ in $G^D$. Assume that a vertex $v$ is contained in $B_2 \cup B_3$. Then, all the two or three white neighbors of $v$ have white-degree $0$ and consequently, they are red vertices in  $G^{D\cup \{v\}}$. It follows that $\sd(v) \ge 3+ 2\cdot 5=13$ that is a contradiction again.

\medskip
(iv) According to (i)-(iii), every non-red vertex of $G^D$ is a blue leaf or a white vertex being adjacent only to blue leaves.  In other words, $G^D$ consists of star-components, each of which contains only one white vertex. If there is a white vertex $v$ which is adjacent to at least three blue leaves, then this component becomes red in $G^{D\cup \{v\}}$ and we have  $\sd(v) \ge 5+ 3\cdot 3=14$. This contradiction establishes (iv).
 \smallqed

\subsection{Phase 4}

Phase~4 starts when Phase~3 ends, but the game is not over. Phase~4 finishes when the game ends.

By Claim~\ref{cl:later2} and \ref{cl:end3}, the residual graph $G^D$ may contain only the following types of components at the beginning of Phase~4: WB$^-$,  WB$^+$, BWB, an isolated red vertex. As follows, playing any (legal) vertex $v$ in the game, exactly one component becomes red. It means $S(v)= 5+3=8$, $S(v)=5+4-1=8$, and $S(v)\ge 5+2\cdot 3 -3=8$ if $v$ belongs to a  WB$^-$-, WB$^+$-, and BWB-component respectively. So, either Staller or Dominator plays a vertex, $F(G^D)$ decreases by at least $8$ with each move.

\begin{claim} \label{cl:av4}
	If Phase~4 consists of $p_4$ moves, then $F(G^D)$ decreases by at least $8p_4$ during this phase. 
\end{claim} 

\subsection{Completion of the proof}

We supposed throughout the proof that Dominator chooses a vertex which results in the maximum achievable decrease in the value of the potential function. For Staller's moves we did not suppose anything but legality.  As the game starts with $f(G^\emptyset)=5n$ and finishes with $F(G^D)=0$,  Claims~\ref{cl:av1}, \ref{cl:av2}, \ref{cl:av3} and \ref{cl:av4} imply that Dominator can ensure
$$5n \ge 8p_1+ 8p_2 + (f(G^{D^{*}})- F(G^{D^{*}}) ) + 8p_3 + 8p_4, $$
where $G^{D^{*}}$ denotes the residual graph obtained at the end of Phase~2 and $p_i$ denotes the number of vertices played in Phase~$i$, for $i=1,\dots 4$. Recall that  $f(G^{D^*})- F(G^{D^*}) \ge 0$. Thus, for the total length $p=p_1+p_2+p_3+p_4$, we have
$$\gamma_g(G) \le p \le \frac{5}{8}\, n.$$
This completes the proof of Theorem~\ref{thm:main}.
\qed

\section{Remark on the Staller-start game}
\label{sec:4}

Already Theorem~\ref{thm:main} and the inequality $|\gamma_g(G) - \gamma_g'(G)| \le 1$ imply that $\gamma_g'(G) \le \frac{5n}{8}+1$ holds for the Staller-start domination game, if $G$ is isolate-free. On the other hand, it is easy to get a bit better estimation by the modification of the proof of Theorem~\ref{thm:main}. If Staller starts the game, we may consider her first move $v_0$ as the $0 ^{\rm th}$ move in the game and declare that it belongs to Phase~1. Then, in the residual graph $G^{\{v_0\}}$, the vertex $v_0$ is red and, since $G$ is isolate-free, at least one further vertex is recolored light blue or red. We conclude $\s(v_0) \ge 5+1=6$. In the continuation, we may follow the line of the proof presented in Section~\ref{sec:3} and get that the $p+1$ moves in the Staller-start domination game fulfill $5n \ge 6+ 8p$. Therefore, every isolate-free graph $G$ satisfies
$$ \gamma_g'(G) \le p+1 \le \frac{5n+2}{8}.$$


\begin{thebibliography}{}
	
	\bibitem{borowiecki-2019}
	M.~Borowiecki, A.~Fiedorowicz, E.~Sidorowicz,
	Connected domination game,
	\emph{Appl.\ Anal.\ Discrete Math.}  13 (2019), 261--289.  
		
	\bibitem{BDKK-2014}
	B. Bre\v{s}ar, P. Dorbec, S. Klav\v{z}ar, and G. Ko\v{s}mrlj,
	Domination game: effect of edge- and vertex-removal.
	\emph{ Discrete Math.} 330 (2014), 1--10.
	
	\bibitem{ZL}
	B. Bre\v{s}ar, Cs.\ Bujt\'as, T. Gologranc S. Klav\v{z}ar, G. Ko\v{s}mrlj, T. Marc, B. Patk\' os, Zs. Tuza, and M. Vizer,
	The variety of domination games.
	\emph{Aequationes Mathematicae} 93 (2019), 1085--1109.
	
	\bibitem{BKKR-2013}
	B.~Bre{\v{s}}ar, S.~Klav{\v{z}}ar, G.~Ko\v smrlj, and D.~F.~Rall,  
	Domination  game:  extremal families  of  graphs  for  the  $3/5$-conjectures.
	\emph{ Discrete Appl. Math.} 161 (2013), 1308--1316.
	
	\bibitem{BKR-2010}
	B.~Bre{\v{s}}ar, S.~Klav{\v{z}}ar, and D.~F.~Rall, Domination game
	and an imagination strategy. \emph{ SIAM J. Discrete Math.} 24 (2010),
	979--991.
	
	\bibitem{bujtas-2015}
	Cs.\ Bujt\'as, Domination game on forests.
	\emph{  Discrete Math.} 338 (2015), 2220--2228.
	
	\bibitem{bujtas-2015-2}
	Cs.~Bujt\'as, On the game domination number of graphs with given
	minimum degree. \emph{ Electron. J. Combin.} 22  (2015), \#P3.29.
	
	\bibitem{bujtas-2018}
	Cs.~Bujt\'as, On the game total domination number.
	\emph{ Graphs Combin.} 34 (2018), 415--425.
	
	\bibitem{BDIK-2019}
	 	Cs.~Bujt\'as, P. Dokyeesun, V.~Ir\v si\v c, and S.~Klav{\v{z}}ar, Connected domination game played on Cartesian products. \emph{Open Math.} 17 (2019), 1269--1280.
	
	\bibitem{BHT-2016}
	Cs.\ Bujt\'as, M. A. Henning, and Zs. Tuza,
	Transversal game on hypergraphs and the $\frac{3}{4}$-conjecture on the total domination game.
	\emph{ SIAM J. Discrete Math.} 30 (2016), 1830--1847.
	
	\bibitem{BHT-2017}
	Cs.\ Bujt\'as, M. A. Henning, and Zs. Tuza,
	Bounds on the game transversal number in hypergraphs.
	\emph{ European J. Combin.} 59 (2017), 34--50.
	
	\bibitem{BPTV}
	Cs.\ Bujt\'as, B. Patk\'os, Zs. Tuza, and M. Vizer,
	Domination game on uniform hypergraphs. 
	\emph{ Discrete Appl. Math.} 258 (2019), 65--75.
	
	\bibitem{BT-2016}
	Cs.\ Bujt\'as and Zs. Tuza, The Disjoint Domination Game. \emph{Discrete Math.} 339 (2016), 1985--1992. 
	
	\bibitem{BT-2019}
	Cs.\ Bujt\'as and Zs. Tuza, Fractional Domination Game.
	\emph{ Electron. J. Combin.} 26  (2019), \#P4.3.
	
	\bibitem{DH-2016}
	P.~Dorbec and  M.~A.~Henning, Game total domination for cycles and paths. \emph{Discrete Appl. Math.} 208 (2016), 7--18.
	
	\bibitem{DKR-2015}
	P. Dorbec, G. Ko\v smrlj, and G. Renault, The domination game played on unions of
	graphs. \emph{  Discrete Math.} 338 (2015), 71--79.
	
	\bibitem{HL-2017}
	M. A. Henning, C. L\" owenstein, Domination game: extremal families for the 3/5-conjecture for forests. 
	\emph{Discuss. Math. Graph Theory} 37 (2017), 369--381.
	
	\bibitem{HK-2016}
	M. A. Henning and W. B. Kinnersley, Domination game: A proof of the 3/5-conjecture
	for graphs with minimum degree at least two. \emph{ SIAM J. Discrete Math.} 30 (2016),
	20--35.
	
	\bibitem{HKR-2015}
	M. A. Henning, S. Klav\v{z}ar, and D. F. Rall,
	Total version of the domination game.
	\emph{ Graphs Combin.}  31 (2015),  1453--1462.
	
	\bibitem{HKR-2017} M. A. Henning, S. Klav\v{z}ar, and D. F. Rall,
	The 4/5 upper bound on the game total domination number.
	\emph{ Combinatorica} 31 (2017), 223--251.
	
	\bibitem{HKR-2018} M. A. Henning, S. Klav\v{z}ar, and D. F. Rall, Game total domination critical graphs. \emph{ Discrete Appl. Math.} 250 (2018), 28--37. 
	
	\bibitem{HR-2016} M. A. Henning and D. F. Rall,
	Progress towards the total domination game
	$\frac{3}{4}$-conjecture.
	\emph{ Discrete Math.} 339   (2016),  2620--2627.
	
	\bibitem{Ir-2019} V. Ir\v si\v c, Effect of predomination and vertex removal on the game total domination number of a graph. \emph{ Discrete Appl. Math.} 257 (2019), 216--225.
	
	\bibitem{JL-2019}
	 Y.~Jiang and M.~Lu, Game total domination for cyclic bipartite graphs. \emph{Discrete Appl. Math.} 265 (2019), 120--127. 
	
	\bibitem{KWZ-2013}
	W.~B.~Kinnersley, D.~B.~West, and R.~Zamani, Extremal problems for
	game domination number. \emph{ SIAM J. Discrete Math.} 27 (2013),
	2090--2107.
	
	\bibitem{KR-2019} S. Klav\v zar and D. F. Rall, Domination game and minimal edge cuts. \emph{ Discrete Math.} 342 (2019), 951--958.
	
	\bibitem{Kos-2017} 
	G. Ko\v smrlj, Domination game on paths and cycles. \emph{Ars Math. Contemp.} 13 (2017), 125--136.
	
	\bibitem{MP-2016}
	N.~Marcus and D.~Peleg,
	The Domination Game: Proving the $3/5$ Conjecture on Isolate-Free Forests.
	arXiv:1603.01181 [cs.DS] (2016).
	
	\bibitem{NSS-2016}
	M. J. Nadjafi-Arani, M. Siggers, and H. Soltani, Characterisation of forests with trivial
	game domination numbers. \emph{ J. Comb. Optim.} 32 (2016), 800--811.
	
	\bibitem{Sch-2016}
	S. Schmidt, The 3/5-conjecture for weakly $S(K_{1,3})$-free forests. \emph{ Discrete Math.} 339
	(2016), 2767--2774.
	
	\bibitem{XL-2018} 	K. Xu and X. Li, On domination game stable graphs and domination game edge-critical graphs. \emph{ Discrete Appl. Math.} 250 (2018), 47--56.
	
	\bibitem{XLK-2018}
	K. Xu, X. Li, and S. Klav\v zar, On graphs with largest possible game domination number.
	\emph{ Discrete Math.} 341 (2018), 1768--1777.
	
	
	
	
	%\bibitem{null} {\komm NOT CITED:}
	%	
	%	T.~W.~Haynes, S.~T.~Hedetniemi, and P.~Slater,
	%	Fundamentals of Domination in Graphs.
	%	Marcel Dekker Inc., New York, NY, 1998.
	
	
\end{thebibliography}
\end{document}